\documentclass[12pt]{article} 
\usepackage{amsfonts,amsmath,amssymb,euscript,amscd}
\usepackage[dvips]{graphics}
 \usepackage{graphicx}
\begin{document}\def\ov{\over} \def\ep{\varepsilon}
\def\inv{^{-1}}\def\be{\begin{equation}} \def\P{\mathbb P}
\def\ee{\end{equation}}\def\x{\xi}  \def\iy{\infty} \def\G{\Gamma}
\def\({\left(} \def\){\right)} \def\la{\lambda} \def\ps{\psi}
\def\ph{\varphi} \def\R{\mathbb R} \def\dl{\delta} \def\t{\tau}
\def\z{\zeta} \def\e{\eta} \def\ch{\raisebox{.3ex}{$\chi$}}
\def\Z{\mathbb Z} \def\pl{\partial} \def\e{\eta} \def\tr{{\rm tr}\,}
\def\l{\ell} \def\O{\Omega} \def\noi{\noindent} \def\Ph{\Phi}
\def\s{\sigma} \def\tn{\otimes} \def\g{\gamma} \def\K{\hat K}
\def\qed{\hfill$\Box$} \def\m{\mu} \def\phy{\ph_\iy} \def\m{\mu} 
\def\KA{K_{\rm Airy}} 
\def\A{{\rm{Ai}}}

\hfill August 9, 2008

\begin{center}{\Large\bf  Asymptotics in ASEP with Step Initial Condition}\end{center}
 
\begin{center}{\large\bf Craig A.~Tracy}\\
{\it Department of Mathematics \\
University of California\\
Davis, CA 95616, USA\\
email: \texttt{tracy@math.ucdavis.edu}}\end{center}

\begin{center}{\large \bf Harold Widom}\\
{\it Department of Mathematics\\
University of California\\
Santa Cruz, CA 95064, USA\\
email: \texttt{widom@ucsc.edu}}\end{center}

\begin{abstract} In previous work the authors considered the asymmetric simple exclusion process on the integer lattice in the case of step initial condition, particles beginning at the positive integers. There it was shown that the probability distribution for the position of an individual particle is given by an integral whose integrand involves a Fredholm determinant. Here we use this formula to obtain three asymptotic results for the positions of these particles. In one an apparently new distribution function arises and in another the distribuion function $F_2$ arises. The latter extends a result of Johansson on TASEP to ASEP.
\end{abstract}

\begin{center}{\bf I. Introduction}\end{center}

In previous work \cite{TW} the authors considered the asymmetric simple exclusion process (ASEP) on the integer lattice $\Z$ in the case of step initial condition, particles beginning at the positive integers $\Z^+$. There it was shown that the probability distribution for the position of an individual particle is given by an integral whose integrand involves a Fredholm determinant. Here we use this formula to obtain three asymptotic results for the positions of these particles.

In ASEP a particle waits exponential time, then moves to the right with probability $p$ if that site is unoccupied (or else stays put) or to the left with probability $q=1-p$ if that site is unoccupied (or else stays put). The formula in \cite{TW} gives the distribution function for $x_m(t)$, the position of the $m$th particle from the left at time~$t$ when all $x_m(0)=m$. 

Here we shall assume that $p<q$, so there is a drift to the left, and establish three results on the position of the $m$th particle when $t\to\iy$. The first gives the asymptotics of the probability $\P(x_m(t)\le x)$ when $m$ and $x$ are fixed; the second, conjectured in \cite{TW}, gives the limiting distribution for fixed $m$ when $x$ goes to infinity; and the third gives the limiting distribution when both $m$ and $x$ go to infinity. In the second result an apparently new distribution function arises and in the third the distribution function $F_2$ of random matrix theory \cite{TW1} arises. (That $F_2$ should arise in ASEP has long been suspected. In the physics literature this is referred to as 
\textit{KPZ universality} \cite{KPZ}.)
   
Before giving the results we state the formula derived in \cite{TW}, valid when $p$ and $q$ are nonzero.
It is given in terms of the Fredholm determinant\footnote{The Fredholm determinant of a kernel $K$ is the operator determinant $\det(I-\la K)$. Properties of these determinants, trace class operators, etc., may be found in \cite{GK}.} of a kernel $K(\x,\,\x')$ on $C_R$, a circle with center zero and large radius $R$ described counterclockwise. It acts as an operator by
\[f(\x)\to\int_{C_R}K(\x,\,\x')\,f(\x')\,d\x',\ \ \ (\x\in C_R).\footnote{All contour integrals are to be given a factor $1/2\pi i$.}\]
We use slightly different notation here, which will simplify formulas later. We set
\[\g=q-p,\ \ \ \t=p/q.\]
The kernel is
\be K(\x,\,\x')=q\,{{\x'}^x\,e^{\ep(\x')t/\g}\ov p+q\x\x'-\x},\label{kernel1}\ee
where
\[\ep(\x)=p\,\x\inv+q\,\x-1.\]
The formula is
\be\P\left(x_m(t/\g)\le x\right)=\int {\det(I-\la K)\ov \prod_{k=0}^{m-1}(1-\la\,\t^k)}\, {d\la\ov \la},\label{P}\ee
where  the integral is taken over a contour enclosing the singularities of the integrand at $\la=0$ and $\la=\t^{-k}\ (k=0,\ldots,m-1)$. We mention here the special case, easily derived from this,
\be\P(x_1(t/\g)>x)=\det(I-K).\label{m=1}\ee

The first formula is concrete.
\vskip1ex
\noi{\bf Theorem 1}. Assume $0<p<q$. For fixed $m$ and fixed $x<m$ we have, as $t\to\iy$, 
\[\P\left(x_m(t/\g)>x\right)\sim \prod_{k=1}^\iy (1-\t^k)\,{t^{2m-x-2}\,e^{-t}\ov (m-1)!\,(m-x-1)!}.\]
 
It is clear probabilistically that $\P\left(x_m(t)>x\right)=0$ for all $t$ when $x\ge m$: for a particle to be to the right of its initial position all particles to its right would have to move simultaneously to the right, which surely has probability zero. This will also be seen in the proof of theorem.

Although Theorem 1 required $p>0$ the statement makes sense when $p=0$, the TASEP where particles move only to the left. In this case the probability equals a probability in a unitary Laguerre random matrix ensemble \cite{Jo}. The corresponding asymptotics can be derived there and found to be the same as our formula when $p=0$.

The second result was conjectured, and the beginning of a possible proof given, in \cite{TW}. Denote by $\K$ the operator on $L^2(\R)$ with kernel\footnote{This is the symmetrization of the Mehler kernel.}
\[\K(z,\,z')={q\ov\sqrt{2\pi}}\,\,e^{-(p^2+q^2)\,(z^2+{z'}^2)/4+pq\,zz'}.\]
 
\noi{\bf Theorem 2}. Assume $0<p<q$. For fixed $m$ the limit
\[\lim_{t\to\iy}\P\({x_m(t/\g)+t\ov\g^{1/2}\,t^{1/2}}\le s\)\]
is equal to the integral in (\ref{P}) with $K$ replaced by the operator $\K\,\ch_{(-s,\;\iy)}$.

From this and (\ref{m=1}) we have the special case
\[\lim_{t\to\iy}\P\({x_1(t/\g)+t\ov\g^{1/2}\,t^{1/2}}>-s\)
=\det\big(I-\K\,\ch_{(s,\;\iy)}\big).\]

This is an apparently new family of distribution functions, parametrized by $p$. When $p=0$ the kernel has rank one and the determinant equals a standard normal distribution.

Finally, we state the result when $m$ and $x$ both go to infinity. 
We use the notations
\be\s=m/t,\ \ \ c_1=-1+2\sqrt{\s},\ \ \ c_2=\s^{-1/6}\,(1-\sqrt{\s})^{2/3}.\label{sigma}\ee

\noi{\bf Theorem 3}. When $0\le p<q$ we have
\[\lim_{t\to\iy}\P\({x_m(t/\g )-c_1\,t\ov c_2\,t^{1/3}}\le s\)=F_2(s)\]
uniformly for $\s$ in a compact subset of $(0,\,1)$.\footnote{Notice that here we allow $p=0$. In this case we get the asymptotic formula derived by Johansson \cite{Jo} for TASEP. For ASEP the strong law $t\inv\,x_m(t/\g)\to c_1$ a.s. was proved by Liggett \cite{Li}. For stationary ASEP Bal\'azs and Sepp\"al\"ainen \cite{BS} and Quastel and Valk\'o \cite{QV} proved that the variance of the current across a characteristic has order $t^{2/3}$ and the diffusivity has order
$t^{1/3}$.}
 
The proofs of the theorems will involve asymptotic analysis of $K$. The main point is that the kernel has the same Fredholm determinant as the sum of two kernels; one has large norm but fixed spectrum and its resolvent can be computed exactly, and the other is better behaved. This representation is derived in the next section.  
\begin{center}{\bf II. Preliminaries}\end{center} 
 
We begin with two facts on stability of the Fredholm determinant. They concern smooth kernels acting on simple closed curves. Both use the fact for a trace class operator $L$ the determinant $\det(I-\la L)$ is determined by the traces $\tr L^n,\ n\in\Z^+$. This is so because up to constants these are the coefficients in the expansion of the logarithm of the determinant around $\la=0$.

\noi{\bf Proposition 1}. Suppose $s\to\G_s$ is a deformation of closed curves and a kernel $L(\e,\,\e')$ is analytic in a neighborhood of $\G_s\times\G_s\subset\mathbb{C}^2$ for each $s$. Then the Fredholm determinant of $L$ acting on $\G_s$ is independent of $s$.

\noi{\bf Proof}. The trace of $L^n$ on $\G_s$ equals
\[\int_{\G_s}\cdots\int_{\G_s}L(\e_1,\,\e_2)\cdots L(\e_{n-1},\,\e_n)\,L(\e_n,\,\e_1)\,d\e_1\cdots d\e_n.\]
If $s'$ is sufficiently close to $s$ we may consecutively replace the contours $\G_s$ for the $\e_i$ by $\G_{s'}$, obtaining the trace of $L^n$ on $\G_{s'}$. So $\tr L^n$ is a locally constant function of $s$ and the usual argument shows that it is constant. Therefore so is the Fredholm determinant.\qed

\noi{\bf Proposition 2}. Suppose $L_1(\e,\,\e')$ and $L_2(\e,\,\e')$ are two kernels acting on a simple closed contour $\G$, that $L_1(\e,\,\e')$ extends analytically to $\e$ inside $\G$ {\bf or} to $\e'$ inside $\G$, and that $L_2(\e,\,\e')$ extends analytically to $\e$ inside $\G$ {\bf and} to $\e'$ inside $\G$. Then the Fredholm determinants of $L_1(\e,\,\e')+L_2(\e,\,\e')$ and $L_1(\e,\,\e')$ are equal.

\noi{\bf Proof}. Suppose $L_1(\e,\,\e')$ extends analytically to $\e'$ inside $\G$. The operator $L_1\,L_2$ on $\G$ has kernel
\[L_1\,L_2\,(\e,\,\e')=\int_\G L_1(\e,\,\z)\,L_2(\z,\,\e')\,d\z=0,\]
since the integrand extends analytically to $\z$ inside $\G$. The operator $L_2^2$ on $\G$ has kernel
\[L_2^2\,(\e,\,\e')=\int_\G L_2(\e,\,\z)\,L_2(\z,\,\e')\,d\z=0\]
for the same reason. Therefore for $n>1$ 
\[(L_1+L_2)^n=L_1^n+L_2^{n-1}\,L_1,\]
and $\tr L_2^{n-1}\,L_1 =\tr L_1\,L_2^{n-1}=0$, so $\tr (L_1+L_2)^n=\tr L_1^n$. When $n=1$ we use 
\[\tr L_2=\int_\G L_2(\e,\,\e)\,d\e=0,\]
since the integrand extends analytically inside $\G$, which completes the proof.\qed

We introduce the notation
\[\ph(\e)=\({1-\t \e\ov1-\e}\)^x\,e^{\left[{1\ov 1-\e}-
{1\ov 1-\t \e}\right]\,t}.\]
In $K(\x,\,\x')$ we make the substitutions
\[\x={1-\t\e\ov1-\e},\ \ \ \x'={1-\t\e'\ov1-\e'},\]
and we obtain the kernel\footnote{This is the kernel $(d\x/d\e)^{1/2}(d\x'/d\e')^{1/2}\,K(\x(\e),\,\x'(\e'))$.} 
\[{\ph(\e')\ov \e'-\t\e}=K_2(\e,\,\e')\]
acting on $\g$, a little circle about $\e=1$ described clockwise, which has the same Fredholm determinant. We denote this by $K_2$ because there is an equally important kernel 
\[{\ph(\t\e)\ov \e'-\t\e}=K_1(\e,\,\e').\]

\noi{\bf Proposition 3}. Let $\G$ be any closed curve going around  $\e=1$ once counterclockwise with $\e=\t\inv$ on the outside.
Then the Fredholm determinant of $K(\x,\,\x')$ acting on $C_R$ has the same Fredholm determinant as $K_1(\e,\,\e')-K_2(\e,\,\e')$ acting on $\G$.

\noi{\bf Proof}. We must show that the determinant of $K_2$ acting on $\g$ equals the determinant of $K_1-K_2$ acting on $\G$. The kernel $K_1(\e,\,\e')$ extends analytically to $\e$ inside $\g$ and to $\e'$ inside $\g$ while $K_2(\e,\,\e')$ extends analytically to $\e$ inside $\g$. Hence by Proposition~2 the determinant of $K_2$ acting on $\g$ equals the determinant of $K_2-K_1$.
Next we show that we may replace $\g$ by $-\G$. (Recall that $\g$ is described clockwise and $\G$ counterclockwise.) We apply Proposition~1 to the kernel 
\[K_1(\e,\,\e')-K_2(\e,\e')={\ph(\t\e)-\ph(\e')\ov \e'-\t\e},\]
with $\G_0=-\g$ and $\G_1=\G$. Since the numerator vanishes when the denominator does, the only singularities of the kernel are at $\e,\,\e'=1,\,\t\inv$, neither of which is passed in a deformation $\G_s,\ s\in[0,\,1]$. Therefore the proposition applies and gives the result.\qed

\noi{\bf Proposition 4}. Suppose the contour $\G$ of Proposition~3 is star-shaped with respect to $\e=0$.\footnote{This means that $0$ is inside $\G$ and each ray from 0 meets $\G$ at exactly one point.} Then the Fredholm determinant of $K_1$ acting on $\G$ is equal to
\[\prod_{k=0}^\iy (1-\la\t^k).\]

\noi{\bf Proof}. The function $\ph(\t\e)$ is analytic except at $\t\inv$ and $\t^{-2}$, both of which are outside $\G$, so the function is analytic on $s\G$ when $0<s\le1$. The denominator $\e'-\t\e$ is nonzero for $\e,\,\e'\in s\G$ for all such $s$. (The assumption on $\G$ was used twice.) Therefore by Proposition~1 the Fredholm determinant of $K_1$ on $\G$ is the same as on $s\G$. This in turn is the same as the Fredholm determinant of
\be s\,K_1(s\e,\,s\e')={\ph(s\t\e)\ov \e'-\t\e}\label{Ks}\ee
on $\G$. The operator is the one with kernel
 \[K_0(\e,\,\e')={1\ov \e'-\t\e},\]
which is trace class since the kernel is smooth, left-multiplied by multiplication by $\ph(s\t\e)$. The latter converges in operator norm to the identity as $s\to0$ since $\ph(s\t\e)\to\nolinebreak1$ uniformly on $\G$, and so (\ref{Ks}) converges in trace norm to $K_0$. Therefore the Fredholm determinant of $K_1$ equals the Fredholm determinant of $K_0$.

The kernel of $K_0^2$ equals
\[K_0^2(\e,\,\e')=\int_\G{d\z\ov (\z-\t \e)\,(\e'-\t \z)}=
{1\ov \e'-\t^2\,\e}.\]
This is because $\t\e$ is inside $\G$ and $\t\inv\e'$ outside $\G$ when when $\e,\,\e'\in\G$, since $\G$ is star-shaped. Generally, we find that
\[K_0^n(\e,\,\e')=\int_\G{d\z\ov (\z-\t^{n-1} \e)\,(\e'-\t \z)}=
{1\ov \e'-\t^n\,\e},\]
so
\[\tr K_0^n={1\ov 1-\t^n}.\]
Thus for small $\la$
\[\log\det(I-\la K_0)=-\sum_{n=1}^\iy{\la^n\ov n}{1\ov1-\t^n}=-\sum_{k=0}^\iy\sum_{n=1}^\iy {\t^{nk}\la^n\ov n}=\sum_{k=0}^\iy\log(1-\la\t^k),\]
and the result follows.\footnote{It is easy to see directly that the nonzero eigenvalues of $K_0$ are exactly the $\t^{-k}$. This does not give the formula for the Fredholm determinant since for that we would have to show that these eigenvalues have algebraic multiplicity one. The computation of traces  avoids that issue.}\qed

Denote by $R(\e,\,\e';\,\la)$ the resolvent kernel of $K_1$, the kernel of $\la\,(I-\la K_1)\inv\,K_1$. This is analytic everywhere except for $\la=\t^{-k},\ k\ge0$. We define
\[\ph_n(\e)=\ph(\e)\,\ph(\t\e)\cdots\ph(\t^{n-1}\e).\]

\noi{\bf Proposition 5}. Assume that $\G$ is star-shaped with $1$ inside and $\t\inv$ outside. Then for sufficiently small $\la$
\[R(\e,\,\e';\,\la)=\sum_{n=1}^\iy \la^n {\ph_n(\t\e)\ov \e'-\t^n\e}.\]

\noi{\bf Proof}. If $0<\t_1,\,\t_2<1$ and $\s_1,\,\s_2$ are analytic inside $\G$ then 
\[\int_\G{\s_1(\e)\ov \z-\t_1\e}\,{\s_2(\z)\ov \e'-\t_2\z}\,d\z=
{\s_1(\e)\,\s_2(\t_1\e)\ov \e'-\t_1\t_2\,\e}.\]
This uses, again, the assumption that $\G$ is star-shaped. From this we see by induction that $K_1^n$ has kernel
\[{\ph_n(\t\e)\ov \e'-\t^n\e}.\]
Here we used the fact that the $\ph_n(\t\e)$ are analytic inside $\G$, although $\ph(\e)$ isn't. We multiply by $\la^n$ and sum to get the resolvent.\qed

For $\la$ not equal to any $\t^{-k}$ the operator  $I-\la K_1$ is invertible and we may factor it out from $I-\la K=I-\la K_1+\la K_2$
and we obtain
\[\det(I-\la K)=\det(I-\la K_1)\,\det\(I+\la K_2\,(I-\la K_1)\inv\)\]
\[=\det(I-\la K_1)\,\det (I+\la K_2\,(I+R))\,,\]
where $R$ denotes the operator with kernel $R(\e,\,\e';\,\la)$.
The first factor is given by Proposition 4 (we asume here that $\G$ is as in the proposition), and so we may rewrite (\ref{P}) as 
\be\P\left(x_m(t/\g)\le x\right)=\int \prod_{k=m}^\iy(1-\la\,\t^k)\cdot \det(I+\la K_2\,(I+R))\, {d\la\ov \la}.\label{P1}\ee
This formula and the formula of Proposition 5 form the basis for our proofs.\footnote{It will become apparent later that it was important that we factored in the order we did. The operator $K_2\,(I+R)$ behaves well while $(I+R)\,K_2$ does not.}

\begin{center}{\bf III. Proof of Theorem 1}\end{center}

We begin this section with a decomposition of the resolvent kernel that will be use in the proofs of the first two theorems. The first summand will contain the poles of the resolvent inside the contour of integration in (\ref{P1}) while the remainder will be analytic inside it. We assume as before that $\G$ is as in Proposition 5.

We have
\[\ph_n(\e)=\({1-\t^n \e\ov1-\e}\)^x\,e^{\left[{1\ov 1-\e}-
{1\ov 1-\t^n \e}\right]\,t},\]
and we define
\[\ph_\iy(\e)=\lim_{n\to\iy}\ph_n(\e)=(1-\e)^{-x}\,e^{{\e\ov 1-\e}t}\]
and
\[G(\e,\,\e',u)=\({1-u \e\ov1-\e}\)^x\,e^{\left[{1\ov 1-\e}-
{1\ov 1-u \e}\right]\,t}\,(\e'-\t\inv u\e)\inv.\]
In this formula we shall always take $u\in [0,\,\t^2]$, so  $G(\e,\,\e',u)$ will be smooth in $u$ and $\e,\,\e'\in\G$. 

Define
\[R_1(\e,\,\e'\,;\la)={\ph_\iy(\t\e)\ov\ph_\iy(\e)}\,\sum_{k=0}^{m-1} {G^{(k)}(\e,\,\e',0)\ov k!}\,{\la\t^{2k}\ov1-\la\t^k},\]
\[R_2(\e,\,\e'\,;\la)={\ph_\iy(\t\e)\ov\ph_\iy(\e)}\sum_{n=1}^\iy \la^n\,{\t^{(n+1)(m-1)}\ov(m-1)!}\int_0^{\t^{n+1}}(1-u/\t^{n+1})^{m-1}\,G^{(m)}(\e,\,\e',u)\,du.\]
(Derivatives of $G$ are all with respect to $u$.) 

Clearly $R_1$ is analytic everywhere except for poles at $\la=1,\,\t\inv,\ldots,\t^{-m+1}$ and $R_2$ is defined and analytic for $|\la|<\t^{-m}$.

\noi{\bf Lemma 1}. $R(\e,\,\e'\,;\la)=R_1(\e,\,\e'\,;\la)+R_2(\e,\,\e'\,;\la)$ when $|\la|<\t^{-m}$.

\noi{\bf Proof}. Observe that 
\[{\ph_n(\t\e)\ov \e'-\t^n\e}={\ph_\iy(\t\e)\ov\ph_\iy(\e)}\,G(\e,\,\e',\t^{n+1}).\]
By Taylor's theorem with integral remainder $G(\e,\,\e',\t^{n+1})$ is equal to
\[\sum_{k=0}^{m-1} {G^{(k)}(\e,\,\e',0)\ov k!}\,\t^{(n+1)k}+
{\t^{(n+1)(m-1)}\ov (m-1)!}\int_0^{\t^{n+1}}(1-u/\t^{n+1})^{m-1}\,G^{(m)}(\e,\,\e',u)\,du.\]
We multiply this by $\ph_\iy(\t\e)/\ph_\iy(\e)$ times $\la^n$ and sum over $n$ to get $R(\e,\,\e'\,;\la)$. We obtain the statement of the proposition for $\la$ sufficiently small, and therefore by analyticity it holds throughout $|\la|<\t^{-m}$\qed

\noi{\bf Lemma 2}. The operators $K_2\,R_1$ and $K_2\,R_2$ have kernels
\be K_2\,R_1(\e,\,\e')=\sum_{k=0}^{m-1}{1\ov k!}\,{\la\t^{2k}\ov1-\la\t^k}\int_\G {G^{(k)}(\z,\,\e',0)\ov\z-\t\e}\,d\z,\label{K2R1}\ee
\be K_2\,R_2(\e,\,\e')=\sum_{n=1}^\iy \la^n\,{\t^{(n+1)(m-1)}\ov (m-1)!}\int_0^{\t^{n+1}}(1-u/\t^{n+1})^{m-1}\,du\int_\G {G^{(m)}(\z,\,\e',u)\ov\z-\t\e}\,d\z.\label{K2R2}\ee

\noi{\bf Proof}. We have
\[\ph(\z)\,{\ph_\iy(\t\z)\ov\ph_\iy(\z)}=1.\]
The formulas (\ref{K2R1}) and (\ref{K2R2}) follow from this and Lemma 1.\qed

When $x$ is fixed the steepest descent curve for $\ph(\e)$ is the  circle with center zero and radius $1/\sqrt{\t}$.
In this section we take for $\G$ the circle with center zero and any radius $r\in (1,\,\t\inv)$, described counterclockwise. This is one of the contours allowed. On $\G$ the function $\ph(\e)$ is well-behaved (it is uniformly exponentially small as $t\to\iy$), but $\ph(\t\e)$ is badly-behaved (it is exponentially large at $\e=r$), which explains the importance of the correct order alluded to in the last footnote.

We begin by deriving trace norm estimates.  
In Lemma 2 the kernels $K_2\,R_1$ and $K_2\,R_2$ are given in terms of integrals of rank one operators, and we shall use the fact that the trace norms of these integrals are at most the integrals of the Hilbert-Schmidt norms of the integrands. We denote by $\|\,\cdot\,\|_1$ the trace norm and (this will be used later) by $\|\,\cdot\,\|_2$ the Hilbert-Schmidt norm. For the estimates involving $R_1$ in the following lemma we assume that $\la$ is bounded away from the poles $\t^{-k}$.

\noi{\bf Lemma 3.} We have, for some $\dl>0$,\footnote{We shall always use $\dl$ to denote some positive number, different with each occurrence.}
\[\|K_2\|_1=O(e^{-\dl t}),\ \ \|K_2\,R_2\|_1=O(e^{-\dl t}),\ \ \|K_2\,R_1\|_1=O(e^{-(1/2+\dl) t}),\]
\be\|K_2\,R_1\,K_2(I+R_2)\|_1=O(e^{-(1+\dl) t}).\label{tracenorms}\ee

\noi{\bf Proof}. For our estimates we use the fact that if $v>0$ then on $\G$ the real part of 
$1/(1-v\e)$
achieves its maximum at $\e=-r$ when $vr>1$ and its minimum at $\e=-r$ when $vr<1$. In particular the real part of
\[{1\ov 1-\e}-{1\ov 1-\t \e}\]
achieves its maximum at $\e=-r$ and equals
\[{1\ov 1+r}-{1\ov 1+\t r}<0.\]

This gives, first, a uniform estimate $\ph(\e)=O(e^{-\dl t})$. The operator $K_2$ equals the operator with trace class kernel $1/(\e'-\t\e)$ left-multiplied by the operator multiplication by $\ph(\e)$, which has operator norm $O(e^{-\dl t})$. This gives the first estimate, $\|K_2\|_1=O(e^{-\dl t})$.

Next, $G^{(m)}(\e,\,\e',u)$ is $O(t^m)$ times the exponential of 
\[\left[{1\ov 1-\e}-{1\ov 1-u \e}\right]\,t,\]
and when $|u|\le\t^{2}$, as it is in (\ref{K2R2}), the real part of this when $\e\in\G$ is at most
\[\left[{1\ov 1+r}-{1\ov 1+\t^2r}\right]\,t,\]
and the expression in brackets is negative. Thus the integrand in the integral over $\z$ in (\ref{K2R2}) is $O(e^{-\dl t})$ uniformly in all variables. In particular its Hilbert-Schmidt norm with respect to $\e,\,\e'$ has the same estimate, so this integral has trace norm $O(e^{-\dl t})$ uniformly in $u$. It follows that $\|K_2\,R_2\|_1=O(e^{-\dl t})$ on compact subsets of $|\la|<\t^{-m}$. 

For $K_2\,R_1$, we use
\be G'(\e,\,\e',\,\,u)=-\left[{\e x\ov1-u\e}+{\e t\ov(1-u\e)^2}-{\t\inv\e\ov \e'-\t\inv u\e}\right]\,G(\e,\,\e',\,u),\label{G'}\ee
from which we see that each
\[{G^{(k)}(\e,\,\e',\,\,0)\ov G(\e,\,\e',\,\,0)}\]
is a linear combination of products $t^i\,\e^j\,(\e')^{-\l}$.
Since
\[G(\e,\,\e',\,\,0)=\ph_\iy(\e)/\e',\]
$G^{(k)}(\e,\,\e',\,\,0)$ is a linear combination of
\[t^i\,\e^j\,\ph_\iy(\e)\,(\e')^{-\l-1},\]
and so by (\ref{K2R1}) $K_2\,R_1(\e,\,\e')$ is a linear combination of integrals
\be t^i\,\int_\G \z^j\,{\ph_\iy(\z)\ov \z-\t\e}\,(\e')^{-\l-1}d\z.\label{K2R1comb}\ee
The exponent in $\ph_\iy(\z)$ is $t$ times
$\z/(1-\z)$.
Its maximum real part on $\G$, occurring at $\z=-r$, is $-r/(1+r)$. Since $r/(1+r)>1/2$ this shows that the integrand is uniformly
$O(e^{-(1/2+\dl) t})$, and so this is the bound for $\|K_2\,R_2\|_1$, as long as $\la$ is bounded away from the poles. 

Finally, $K_2\,R_1\,K_2$ and $K_2\,R_1\,K_2\,R_2$. It follows from (\ref{K2R1comb}) that the kernel of $K_2\,R_1\,K_2$ is a linear combination of
\[t^i\int_\G\int_\G \z^j\,{\ph_\iy(\z)\ov \z-\t\e}\,(\z')^{-\l-1}\,{\ph(\e')\ov \e'-\t\z'}\,d\z\,d\z'.\]
We integrate first with respect to $\z'$ by expanding the contour. We cross the pole at $\z'=\t\inv \e'$, so we get a constant times
\[t^i\int_\G {\ph_\iy(\z)\,\ph(\e')\ov\z-\t\e}\,\z^j\,(\e')^{-\l-1}d\z.\]
Now we compute
\[\int_\G{\ph_\iy(\z)\ov\z-\t\e}\,\z^j\,d\z.\]
The integrand is analytic outside $\G$ with a pole at infinity. The integral may be written
\[\sum_{k=0}^\iy (\t\e)^k\,\int_\G (1-\z)^{-x}\,e^{{\z\ov 1-\z}t}\,\z^{j-k-1}\,d\z,\]
which we see equals $e^{-t}$ times a polynomial in $t$ and 
$\e$. So the kernel of $K_2\,R_1\,K_2$ is $e^{-t}$ times a linear combination of products $t^i\,\e^j\,\ph(\e')\,(\e')^{-\l-1}$. 

Since $\ph(\e')=O(e^{-\dl t})$ as we have already seen, we have $\|K_2\,R_1\,K_2\|_1=O(e^{-(1+\dl) t})$.
If we use $\ph(\e)\,\ph_\iy(\t\e)=\ph_\iy(\e)$ again we see that the kernel of $K_2\,R_1\,K_2\,R_2$ is $e^{-t}$ times a linear combination of 
\[t^i\, \e ^j\,\sum_{n=1}^\iy \la^n\,{\t^{(n+1)(m-1)}\ov(m-1)!}\int_\G\int_0^{\t^{n+1}}\z^{-\l-1}\,(1-u/\t^{n+1})^{m-1}\,G^{(m)}(\z,\,\e',u)\,du\,d\z.\]
Using (\ref{G'}) again we see that $G^{(m)}(\z,\,\e',u)$ is $O(t^m)$ times the exponential of
\[\left[{1\ov 1-\z}-{1\ov 1-\z u}\right]\,t.\]
As before the maximum real part of the expression in brackets occurs at $-r$ and equals
\[{1\ov 1+r}-{1\ov 1+r u},\]
which has a negative upper bound for $u\le\t^2$. Since we had the factor $e^{-t}$ we obtain the bound $\|K_2\,R_1\,K_2\,R_2\|_1=O(e^{-(1+\dl) t})$.

This completes the proof of Lemma 3.\qed

\noi{\bf Proof of Theorem 1}. In (\ref{P1}) the contour encloses  all the singularities of the integrand. If we take the contour instead to have the singularity $\la=0$ on the outside and the $\t^{-k}$ with $k<m$ inside then we have
\be\P\left(x_m(t/\g)>x\right)=-\int \prod_{k=m}^\iy(1-\la\,\t^k)\cdot \det(I+\la K_2\,(I+R))\, {d\la\ov \la}.\label{P3}\ee

Now
\[I+\la K_2\,(I+R)=I+\la K_2\,(I+R_2)+\la K_2\,R_1\]
\[=
(I+\la\,K_2\,R_1\,(I+\la\,K_2(1+R_2))\inv)\,(I+\la\,K_2(1+R_2)).\]
(Note that $I+\la\,K_2(1+R_2)$ is invertible since $K_2(1+R_2)$ has small norm.) Therefore
\[\det(I+\la K_2\,(I+R))=\det(I+\la\,K_2(1+R_2))\,\det(I+\la\,K_2\,R_1\,(I+\la\,K_2(1+R_2))\inv).\]
The first factor on the right is analytic inside the contour, and equal to $1+O(e^{-\dl t})$ by (\ref{tracenorms}). As for the second factor, we have
\[I+\la\,K_2\,R_1\,(I+\la\,K_2(1+R_2))\inv=I+\la\,K_2\,R_1-\la^2 K_2\,R_1\,K_2(1+R_2)(I+\la\,K_2(1+R_2))\inv\]
\[=I+\la\,K_2\,R_1+O(e^{-(1+\dl)t}),\]
by (\ref{tracenorms}). Here the error estimate refers to the trace norm. Hence 
\[\det(I+\la\,K_2\,R_1\,(I+\la\,K_2(1+R_2))\inv)=1+\la\,\tr K_2\,R_1+O(e^{-(1+\dl)t}),\]
since $\|K_2\,R_1\|_1=O(e^{-(1/2+\dl)t})$ by (\ref{tracenorms}), so $\|(K_2\,R_1)^2\|_1=O(e^{-(1+\dl)t})$. Thus
\[\det(I+\la K_2\,(I-\la K_1)\inv)=\det(I+\la\,K_2(1+R_2))\,\(1+\la\,\tr K_2\,R_1+O(e^{-(1+\dl)t})\).\]

When we insert this into the integral in (\ref{P3}) we may ignore the summand $1$ in the second factor since the first factor is analytic inside the contour. The integral involving $\tr K_2\,R_1$ we can compute by residues. Its multiplier $\la$ is cancelled by the denominator in (\ref{P3}). So with error $O(e^{-(1+\dl)t})$ (\ref{P3}) equals
\be-\sum_{k=1}^{m-1}\,\prod_{j=m}^\iy (1-\t^{j-k})\cdot
\det(I+\la\,K_2(1+R_2(\t^{-k}))\cdot
\textrm{residue of}\ \tr K_2\,R_1\ {\rm at}\ \la=\t^{-k},\label{approx}\ee
where $R_2(\t^{-k})$ denotes the operator with kernel $R_2(\e,\,\e';\,\t^{-k})$.

The determinants are $1+O(e^{-\dl t})$, as we saw, and will not contribute to the asymptotics. The residue of $\tr K_2\,R_1$ at $\t^{-k}$ equals
\be-{1\ov k!}\int_\G\int_\G{G^{(k)}(\z,\,\e,0)\ov\z-\t\e}\,d\z\,d\e.\label{kres}\ee
{}From (\ref{G'}) we see, more precisely than earlier, that 
\[ G^{(k)}(\z,\,\e,0)=\ph_\iy(\z)\,\z^{k}\,\sum_{i+j\le k}a_{ijk}\,t^i\,\e^{-j-1}\]
for some coefficients $a_{ijk}$.
Substituting this into (\ref{kres}) and integrating with respect to $\e$ by expanding the contour outward gives
\[-{1\ov k!}\sum_{i+j\le k} a_{ijk}\,t^i\,\t^{j}\int_\G \ph_\iy(\z)\,\z^{k-j-1}\,d\z
=-{1\ov k!}\sum_{i+j\le k} a_{ijk}\,t^i\,\t^j\int_\G(1-\z)^{-x}\,e^{{\z\ov 1-\z}t}\,\z^{k-j-1}\,d\z.\]
The integral vanishes unless $x\le k-j$ and otherwise equals $e^{-t}$ times a polynomial in $t$ of degree $k-j-x$ with top coefficient
\[{(-1)^{j-k}\ov(k-j-x)!}.\]
We see from this that the highest power of $t$, which is $t^{2k-x}$, comes from the summand with $j=0,\,i=k$. The coefficient $a_{k,0,k}$ equals $(-1)^k$. Thus (\ref{kres}) equals $e^{-t}$ times a polynomial of degree $2k-x$ in $t$ with top coefficient 
\[-{1\ov k!\,(k-x)!}.\]
In particular the main contribution to the sum in (\ref{approx}) comes from the summand $k=m-1$, and if we recall the minus sign in (\ref{approx}) we get the statement of Theorem~1.\qed

\noi{\bf Remark}. As mentioned in the introduction, we can also show that
$\P(x_m(t)>x)=0$ when $x\ge m$. We know for (\ref{K2R1comb}) that $K_2\,R_1(\e,\,\e')$ is a linear combination of
\[(\e')^{-j-1}\int_\G{\ph_\iy(\z)\ov \z-\t\e}\,\z^k\,d\z=
(\e')^{-j-1}\int_\G{1\ov \z-\t\e}\,(1-\z)^{-x}\,e^{{\z\ov 1-\z}t}\,\z^{k}\,d\z,\]
with $j,\,k<m$. When $x\ge m$ we expand the contour and get zero since $k<m$ and $\t\e$ is inside $\G$. Therefore $K_2\,R_1=0$. Hence
\[K_2\,(I-\la K_1)\inv=K_2(1+R_2)+K_2\,R_1=K_2(1+R_2),\]
and so
\[\det\(I+\la K_2\,(I-\la K_1)\inv\)=\det\,(I+\la\,K_2(1+R_2)),\] 
which is analytic inside the contour of integration and therefore integrates to zero.
\pagebreak

\begin{center}{\bf IV. Proof of Theorem 2}\end{center}

We know from (\ref{P}) that
\[\P\({x_m(t/\g)+t\ov\g^{1/2}\,t^{1/2}}\le s\)=\int {\det(I-\la K)\ov \prod_{k=0}^{m-1}(1-\la\,\t^k)}\, {d\la\ov \la},\]
where in the definition of $K$ we set
\be x=-t+\g^{1/2}\,s\,t^{1/2}.\label{x}\ee
Therefore the theorem would follow if 
\be \lim_{t\to\iy}\det(I-\la K)=\det\big(I-\la\K\,\ch_{(-s,\,\iy)}\big)\label{detlim}\ee
uniformly on compact $\la$-sets.
The Fredholm determinants are entire functions of $\la$, and the coefficients in their expansions about $\la=0$ are universal polynomials in the traces of powers of the operators. It was shown in \cite{TW} that for $n\in\Z^+$
\[\lim_{t\to\iy}\tr K^n=\tr \big(\K\,\ch_{(-y,\,\iy)})\big)^n,\] 
and it was pointed out that (\ref{detlim}) would follow if we knew that $\det(I-\la K)$ is uniformly bounded for large $t$ on compact $\la$-sets. This is what we shall show here.

For any $m$ it suffices that the determinant is uniformly bounded on compact subsets of $|\la|<\t^{-m}$, and since it is entire we may assume that the sets exclude the singularities at $\la=\t^{-k}$. From the uniform boundedness of $\det(I-\la K_1)$ on compact $\la$-sets it follows that it suffices to prove the uniform boundedness of $\det(I+K_2\,(I+R))$ on on compact sets excluding the $\t^{-k}$.

Here is how we decide what contour to take for $\G$. The steepest descent curves for all the $\ph_n$ including $\ph_\iy$ are similar. They lie in the right half-plane, tangent to the imaginary axis at the saddle point $\e=0$, and have an inward-pointing cusp at $\e=1$ where the real part of the exponential tends to $-\iy$. 
We would like to take as the curve $\G$ of Propositions 3 and 4 something like this. It need not have that cusp at $\e=1$, only that the $\ph_n$ are exponentially small there, and if it passes through $\e=1$ vertically that will happen. That $\e'=1$ is a singularity of $K_2$ does not change the conclusions of the propositions since we can take an appropriate limit of contours not passing through 1. So we may take $\G$ to be the circle with diameter $[0,\,1]$. But this is not star-shaped with respect to the origin, so Proposition 4 would not apply (even though Proposition 3 would). Therefore we expand it a little on the left, resulting in a contour that is star-shaped. We expand it so that instead of 0 it passes through $-t^{-1/2}$.  This, finally, is the contour $\G$ in  this section: the circle symmetric about the real line and meeting it at $\e=-t^{-1/2}$ and $\e=1$.

From the identity
\[\det(I+A)={\det}_2(I+A)\,e^{\tr A}\] 
and the fact that the ${\det}_2$ is bounded on $\|\,\cdot\,\|_2$-bounded sets, we see that is suffices to prove that 
\[\tr\(K_2\,(I+R)\)=O(1),\ \ \ \|K_2\,(I+R)\|_2=O(1).\]
We shall prove more, namely

\be \tr K_2=O(1),\ \ \ \|K_2\|_2=O(1),\ \ \ \|K_2\,R\|_1=O(1).\label{KRestimates}\ee

We begin by obtaining a bound for integrals involving the various $\ph_n(\e)$. The coefficients of $t$ appearing in the exponentials of these functions are of the form
\be{1\ov1-\e}-{1\ov1-v\e}+\log{1-\e\ov1-v\e}\label{form}\ee
with $0\le v\le\t$. On the part of $\G$ outside any fixed neighborhood of zero in $\mathbb{C}$ the real parts of these are uniformly bounded above by $-\dl$ for some $\dl>0$ when $t$ is sufficiently large.In a sufficiently small fixed neighborhood of zero the real part is at most $O(t\inv)-\dl\,|\e|^2$. It follows that  $\ph_n(\e)=O(e^{-\dl|\e|^2t+O(t^{1/2}|\e|)})$, where the $t^{1/2}|\e|$ term comes from the $y\,t^{1/2}$ term in (\ref{x}). From this it follows that for any $k\ge0$
\be\int_\G|\ph_n(\e)|\,|\e|^k\,|d\e|=O(t^{-(k+1)/2}),\label{phint}\ee
for the following reason. The integral over that part of $\G$ outside any fixed neighborhood of zero is exponentially small. For the integral over a neighborhood of zero we have, if $y={\rm Im}\,\e$,
\[|\e|^2=O(t\inv+y^2),\ \ \ |\e|^2\ge y^2,\ \ \ |d\e|=O(dy),\]
so the integral over that portion of $\G$ is bounded by a constant times
\[\int_{-\iy}^\iy e^{-\dl y^2\,t+O(|y|\,t^{1/2})}\,(t\inv+y^2)^{k/2}\,dy=O(t^{-(k+1)/2}).\]  

If we change variables in (\ref{phint}) we get the equivalent estimate
\[\int_{t^{1/2}\G}|\ph_n(t^{-1/2}\e)|\,|\e|^k\,|d\e|=O(1).\]
More generally, for all $j>0$ we have
\be\int_{t^{1/2}\G}|\ph_n(t^{-1/2}\e)|^j\,|\e|^k\,|d\e|=O(1),\label{estimate}\ee
since $\ph_n(\e)$ is uniformly bounded.

We shall now establish (\ref{KRestimates}). First $K_2$, with kernel
\[{\ph(\e')\ov\e'-\t\e}.\]
We use the fact that the kernel substitution
\be L(\e,\,\e')\ {\rm on}\ \G\ \longrightarrow \ t^{-1/2}\,L(t^{-1/2}\e,\,t^{-1/2}\e')\ {\rm on}\ t^{1/2}\G\label{kernsubs}\ee
preserves norms and traces. The circle $t^{1/2}\G$ meets the real line at $\e=-1$ and $\e=t^{1/2}$. Making this substitution gives the kernel 
\be{\ph(t^{-1/2}\e')\ov\e'-\t\e}.\label{newK2}\ee

We have 
\[|\tr K_2|\le{1\ov1-\t}\int_{t^{1/2}\G}\left|{\ph(t^{-1/2}\e)\ov\e}\right|\,|d\e|=O(1),\]
by (\ref{estimate}) and the fact that $t^{1/2}\G$ is bounded away from zero.

Next,
\[\|K_2\|_2^2=\int_{t^{1/2}\G}\int_{t^{1/2}\G}\left|{\ph(t^{-1/2}\e')\ov\e'-\t\e}\right|^2\,d\e\,d\e'.\]
Now
\[\int_{t^{1/2}\G}{1\ov |\e'-\t\e|^2}\,|d\e|=O(1)\]
uniformly for $\e'\in t^{1/2}\G$.\footnote{That's because if $\e\in  t^{1/2}\G$ then the distance from $\t\e$ to $t^{1/2}\G$ is at least some positive constant times $|\e|$.} Using this and (\ref{estimate}) we see that $\|K_2\|_2=O(1)$. 

Next, $K_2\,R$. 

When $x$ is given by (\ref{x}) we find that
\[G'(\e,\,\e',\,\,u)=-\left[{u\e^2t\ov(1-u\e)^2}+{(q-p)^{-1/2}y\e t^{1/2}\ov 1-u\e}-{\t\inv\e\ov \e'-\t\inv u\e}\right]\,G(\e,\,\e',\,u).\]
{}From this and the fact that $u\e$ is bounded away from 1 when $\e\in\G$ and $u\le\t^2$ we find that each
\[{G^{(k)}(\e,\,\e',\,\,u)\ov G(\e,\,\e',\,\,u)}\]
is bounded by a linear combination of products
\[\left|\e\,t^{1/2}\right|^i\,\left|{\e\ov \e'-u\e/\t}\right|^j.\]

Since $G(\e,\,\e',\,0)=\ph_\iy(\e)/\e'$ it follows in particular that $G^{(k)}(\e,\,\e',\,0)$ is bounded by a constant times a linear combination of products
\[|\e\,t^{1/2}|^i\,|\e|^j\,|\e'|^{-j-1}\,|\ph_\iy(\e)|.\]
After the substitution (\ref{kernsubs}) and the variable change $\z\to t^{-1/2}\z$ in each integral in (\ref{K2R1}) we get as bound a linear combination of
\be\int_{t^{1/2}\G}\left|{\ph_\iy(t^{-1/2}\z)\ov\z-\t\e}\right|\,|\z|^{i+j}|\,\e'|^{-j-1}|\,|d\z|.\label{traceest}\ee
The Hilbert-Schmidt norm with respect to $\e,\,\e'$ of $|\z-\t\e|\inv\,|\e'|^{-j-1}$ on $t^{1/2}\G$ is uniformly bounded for $\z\in t^{1/2}\G$ (as in the last footnote), and so the trace norm of the integral is bounded by
\[\int_{t^{1/2}\G}|\ph_\iy(t^{-1/2}\z)|\,|\z|^{i+j}|\,|d\z|=O(1),\]
by (\ref{estimate}).

That takes care of $K_2\,R_1$. For $K_2\,R_2$ it is enough to show that the last integral in (\ref{K2R2}) has bounded trace norm for $u\le\t^2$, for then the trace norm of $K_2\,R_2$ would be at most a constant times $\sum_{n=1}^\iy |\t^m\la|^n,$ which is bounded on compact subsets of $|\la|<\t^{-m}$.

In the estimate for the last integral the analogue of (\ref{traceest}) would be
\[\int_{t^{1/2}\G}\left|{G_0(t^{-1/2}\z,u)\ov\z-\t\e}\right|\,|\z|^{i+j}\,|\e'-u\z/\t|^{-j-1}|\,|d\z|,\]
where 
\[G_0(\e,\,u)=\({1-u \e\ov1-\e}\)^x\,e^{\left[{1\ov 1-\e}-
{1\ov 1-u \e}\right]\,t}.\]
(This is $G$ without its last factor.) 
Taking the Hilbert-Schmidt norm with respect to $\e,\,\e'$ under the integral sign shows (as in the last footnote again) that the trace norm of the integral is bounded by
\[\int_{t^{1/2}\G}|G_0(t^{-1/2}\z,\,u)|\,|\z|^{i+j}|\,|d\z|.\]
In $G_0$ the factor of $t$ in the exponent is of the form (\ref{form}) with $v=u$, and so this integral is $O(1)$ uniformly for $u\le\t^2$. 
 
This completes the proof of (\ref{KRestimates}) and so of Theorem 2.\qed

\begin{center}{\bf V. Proof of Theorem 3}\end{center} 

In formula (\ref{P1}) the integral is taken over a circle with center zero and radius larger than $\t^{-m+1}$. We set 
\be\la=\t^{-m}\,\m,\label{lamu}\ee
and the formula becomes
\be\P(x_m(t/\g)\le x)=\int \prod_{k=0}^{\iy}(1-\m\,\t^k)\,\cdot\,\det \(I+\t^{-m}\,\m\, K_2\,(I+R)\) \;{d\m\ov \m},\label{P2}\ee
where $\m$ runs over a circle of fixed radius larger than $\t$ (but not equal to any $\t^{-k}$ with $k\ge0$). We shall show that when $c_1$ and $c_2$ are given by (\ref{sigma}) and
\be x=c_1\,t+c_2\,s\,t^{1/3}\label{x1}\ee  
the determinant in this integrand has the limit $F_2(s)$  uniformly in $\m$ and $\s$, which will establish the theorem.  

The main lemma replaces the kernel $\t^{-m}\,\m\, K_2\,(I+R)$ by one which will allow us to do a steepest descent analysis. Now we do not decompose $R$ into a sum of two kernels, but use the entire infinite series in Proposition~5.

We define
\[f(\m,\,z)=\sum_{k=-\iy}^\iy{\t^{k}\ov 1-\t^{k}\m}\,z^k.\]
This is analytic for $1<|z|<\t\inv$ and extends analytically to all $z\ne0$ except for poles at the $\t^k,\ k\in\Z$.
We define a kernel $J(\e,\,\e')$ acting on a circle with center zero and radius $r\in(0,\,1)$ by
\be J(\e,\,\e')=\int{\phy(\z)\ov\phy(\e')}\,{\z^m\ov (\e')^{m+1}}\,{f(\m,\z/\e')\ov\z-\e}\;d\z,\label{J}\ee
where the integral is taken over a circle with center zero and radius in the interval $(1,\,r/\t)$.

\noi{\bf Lemma 4}. With $\la$ given by (\ref{lamu}) we have
\[\det\(I+\la\, K_2\,(I+R)\)=\det(I+\mu\,J).\]

\noi{\bf Proof}. Our operators $K_1$ and $K_2$ may be taken to act on a circle with radius $r\in (1,\,\t\inv)$. From Proposition~5 and the identity
\[\ph_n(\z)={\phy(\z)\ov\phy(\t^n\z)}\]
we obtain
\[K_2\,R(\e,\,\e')=\sum_{n=1}^\iy\la^n\,\int{\phy(\z)\ov
\phy(\t^{n+1}\z)}\,{d\z\ov(\z-\t\e)\,(\e'-\t^{n}\z)}.\]
Here $|\z|=r$ but by analyticity we may take any radius such that
\[1<|\z|<\t^{-1}r.\]
This is equal to 
\[\sum_{n=1}^\iy\la^n\,\int{\phy(\z)\ov\z-\t\e}\,d\z\,
\int{1\ov\phy(u\z)\,(\e'-u\z/\t)}\;{du\ov u-\t^{n+1}},\]
as long as on the circle of $u$-integration we have
\[\t^2<|u|<\t r/|\z|.\]
We use   
\[{1\ov u-\t^{n+1}}=\sum_{k=0}^\iy{\t^{(n+1)k}\ov u^{k+1}}\] 
and sum over $n$ first to get
\[\sum_{k=0}^\iy{\t^{2k}\la\ov 1-\t^k\la}
\int{\phy(\z)\ov\z-\t\e}\,d\z\,\int{1\ov\phy(u\z)\,(\e'-u\z/\t)}\;{du\ov u^{k+1}}.\]

If we assume also that
\be\t<|u|<\t r/|\z|,\  \textrm{which requires also that}\ 1<|\z|<r,\label{domain}\ee
we may rewrite this as
\[\sum_{k=0}^\iy{\t^{k}\la\ov 1-\t^k\la}
\int{\phy(\z)\ov\z-\t\e}\,d\z\,\int{1\ov\phy(u\z)\,(\e'-u\z/\t)}\;{du\ov u^{k+1}}\]
\[-\sum_{k=0}^\iy \t^{k}
\int{\phy(\z)\ov\z-\t\e}\,d\z\,\int{1\ov\phy(u\z)\,(\e'-u\z/\t)}\;{du\ov u^{k+1}},\]
because both series converge.

Summing the second series gives  
\[-\int{\phy(\z)\ov\z-\t\e}\,d\z\,\int{du\ov\phy(u\z)\,(\e'-u\z/\t)\,(u-\t)}.\]
Since $\phy(u\z)$ is analytic and nonzero inside the $u$-contour (since $|u\z|<\t r<1$) and $\t$ is inside and $\t\e'/\z$ outside this equals
\[-\int{\phy(\z)\ov\phy(\t\z)}\,{d\z\ov(\z-\t\e)\,(\e'-\z)}
=-\int{\ph(\z)\ov (\z-\t\e)\,(\e'-\z)}\,d\z.\]
If we expand the contour so that
\[|\z|>r\]
then we pass the pole at $\z=\e'$ and get
\[-{\ph(\e')\ov \e'-\t\e}-\int_{|\z|>r}{\ph(\z)\ov (\z-\t\e)\,(\e'-\z)}\,d\z.\]

The first summand is exactly $-K_2(\e,\,\e')$, so have shown
\[ K_2(I+R)\,(\e,\,\e')=-\int_{|\z|>r}{\ph(\z)\ov (\z-\t\e)\,(\e'-\z)}\,d\z\]
\be+\sum_{k=0}^\iy{\t^{k}\ov 1-\t^k\la}
\int{\phy(\z)\ov\z-\t\e}\,d\z\,\int{1\ov\phy(u\z)\,(\e'-u\z/\t)}\;{du\ov u^{k+1}}.\label{op1}\ee

If the index $k$ were negative then the $u$-integration would give zero since the integrand would be analytic inside the $u$-contour. Therefore the sum over $k$ can be taken from $-\iy$ to $\iy$.

The integration domains in the double integral are given in (\ref{domain}).
If we make the variable change $u\to u/\z$ in the integral the sum becomes 
\[ J_0(\e,\,\e')=\sum_{k=-\iy}^\iy{\t^{k}\ov 1-\t^k\la}
\int{\phy(\z)\ov\z-\t\e}\,\z^k\,d\z\,\int{1\ov\phy(u)\,(\e'-u/\t)}\;{du\ov u^{k+1}},\]
and the new conditions are  
\[1<|\z|<r,\ \ \ \t|\z|<|u|<\t r.\] 

The first operator on the right side of (\ref{op1}) is analytic for $|\e|,\,|\e'|\le r$. The kernel $J_0(\e,\,\e')$ is analytic for $|\e|\le r$. It follows by Proposition~2 that the Fredholm determinant of the sum of the two, i.e., of $K_2(I+R)$, equals the Fredholm determinant of $J_0$.

Now we use (\ref{lamu}). Substituting $k\to m+k$ in the first sum below we find that 
\[\sum_{k=-\iy}^\iy{\t^{k}\ov 1-\t^k\la}\({\z\ov u}\)^k=\t^m\,\({\z\ov u}\)^m\,
\sum_{k=-\iy}^\iy{\t^{k}\ov 1-\t^{k}\m}\({\z\ov u}\)^k=\t^m\,\({\z\ov u}\)^m\,f(\m,\z/u).\]
Thus 
\[J_0(\e,\,\e')=\t^m\,\int\int{\phy(\z)\ov\phy(u)}\,\({\z\ov u}\)^m\,{f(\m,\z/u)\ov(\z-\t\e)\,(\e'-u/\t)}\;d\z\;{du\ov u}.\]

This has the same Fredholm determinant as
\[\t\inv\,J_0(\t\inv\e,\,\t\inv\e')=\t^m\,\int\int{\phy(\z)\ov\phy(u)}\,\({\z\ov u}\)^m\,{f(\m,\z/u)\ov(\z-\e)\,(\e'-u)}\;d\z\;{du\ov u},\]
where now the operator acts on a circle with radius $r\in (\t,\,1)$ and in the integral
\[1<|\z|<r/\t,\ \ \ \t|\z|<|u|<r.\]

We now do something similar to what we did before. If we move the $u$-integral outward, so that $r<|u|<1$ on the new contour, we pass the pole at $u=\e'$, which gives the contribution
\[\t^m\,\int{\phy(\z)\ov\phy(\e')}\,{\z^m\ov (\e')^{m+1}}\,{f(\m,\z/\e')\ov\z-\e}\;d\z=\la\inv\m\, J(\e,\,\e').\]
(The function $f(\m,\,\z/u)$ remains analytic in $u$ during the deformation.) The new double integral is a kernel analytic for $|\e|,\,|\e'|\le r$ and $J(\e,\,\e')$ is analytic for $|\e|\le r$. Therefore by Proposition~2
\[\det\(I+\la\, K_2\,(I+R)\)=\det(I+\mu\,J).\]
The conditions on $r$ and $|\z|$ under which we obtained this were $r\in(\t,\,1)$ and 
\linebreak
$1<|\z|<r/\t$. These may be relaxed as in the statement of the lemma because we may deform the $\e$- and $\z$-contours without passing a singularity of the integrand, and use Proposition~1.\qed 

\noi{\bf Remark}. The lemma was proved under the assumption that $\t>0$. The only occurrence of $\t$ in $\m\,J(\e,\,\e')$ is in $\m\,f(\z/\e')$ and as $\t\to0$ this tends to $\z/(\e'-\z)$. Since the probabilities $\P(x_m(t)\le x)$ are continuous in $p$ at $p=0$\footnote{This follows, for example, from formula (2) of \cite{TW}.} the integral fomula we derived for the probability holds for $p=0$ as well, with this replacement for $\m\,f(\z/\e')$. The asymptotics that follow are actually simpler in this case.
 
We now explain where the constants $c_1$ and $c_2$ come from.
When we to do a saddle point analysis of the integral in (\ref{J}) the first step is to write $\phy(\z)\,\z^m$ as the exponential of
\[-x\,\log(1-\z)+t\,{\z\ov1-\z}+m\,\log\z,\]
and differentiate this to get the saddle point equation
\[{x\ov1-\z}+{t\ov(1-\z)^2}+{m\ov\z}=0,\]
or
\[(m-x)\,\z^2+(x+t-2m)\,\z+m=0.\]
The transition of the asymptotics occurs when the two saddle points coincide, which is when
\[(x+t-2m)^2=4\,m\,(m-x).\]
This gives
\[m={(x+t)^2\ov4t}.\]
Setting $m=\s t$ and $x=c_1\,t$ gives
\[\s={(c_1+1)^2\ov4},\]
or $c_1=-1\pm 2\sqrt{\s}$.
Since $c_1$ should be increasing with $\s$ we take the positive square root in (\ref{sigma}). The saddle point is at
\[\x=-\sqrt\s/(1-\sqrt{\s}).\]
 We compute that if $x$ is given by 
(\ref{x1}) precisely and we set
\[\phy(\z)\,\z^m=\phy(\x)\,\x^m\,e^{\ps(\z)},\]
then in a neighborhood of $\z=\x$
\be\ps(\z)=-c_3^3\,t\,(\z-\x)^3 /3+c_3\,s\,t^{1/3}\,(\z-\x)+O(t(\z-\x)^4))+O(t^{1/3}\,(\z-\x)^2),\label{psi}\ee
where
\[c_3=\s^{-1/6}\,(1-\s^{1/2})^{5/3}.\]
(It is only with $c_2$ as given in (\ref{sigma}) that the coefficients of $t$ and $t^{1/3}$ are related this way.)

Carrying out the details, we define 
\[\ps_0(\z)=-c_1\,\log(1-\z)+{\z\ov1-\z}+\s\log\z,\ \ \ \ps_1(\z)=\ps_0(\z)-\ps_0(\x).\]
There are two steepest descent curves, an outer one $\G_o$ and an inner one $\G_i$. (See Fig.~\ref{sdCurve}. All curves are for the case $\s=1/4$.) Both pass through $\x$ and have cusps at 1. The outer one emanates from $\x$ in the directions $\pm 2\pi/3$ and has an inner-pointing cusp at $\z=1$. On it, Re$(\ps_1(\z))$ has its maximum of zero at $\z=\x$ and tends to $-\iy$ at the cusp. The inner one emanates from $\x$ in the directions $\pm \pi/3$ and has an outer-pointing cusp at $\e=1$. On it, Re$(\ps_1(\e))$ has its minimum of zero at $\e=\x$ and tends to $+\iy$ at the cusp.
\addtocounter{figure}{-1}
\begin{figure} 
\begin{center}
 \includegraphics[scale=.75]{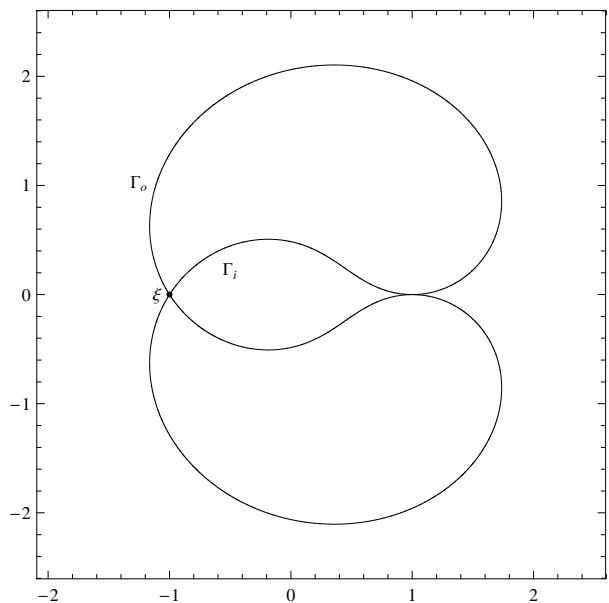}
 \caption{Steepest descent curves $\Gamma_o$ and $\Gamma_i$ for $\psi_1$. The point $\xi$ is the location of the saddle point.\label{sdCurve}}
 \end{center}
 \end{figure}

We would like to deform the $\e$-contour for $J$, which is a circle with radius $r<1$, to $\G_i$ and apply Proposition~1 to assure that the Fredholm determinant doesn't change. The $\z$-contour started out as a circle with radius slightly bigger than one. We may deform the $\e$-contour as described if we deform the $\z$-contour simultaneously, assuring that the $\z$-contour is always just outside the $\e,\,\e'$-contour, so that in particular we don't pass a singularity of $f(\m,\,\z/\e')$. Next we want to expand the $\z$-contour outward to $\G_o$, but in the process we might encounter a singularity of $f(\m,\,\z/\e')$, and this causes a problem. It will happen if a ray from zero meets $\G_i$ at a point $\e$ and $\G_0$ at $\z$ and $\e/\z\le \t$. This will not happen if $\t$ is close enough to zero but will happen if $\t$ is close enough to one.

But we do not have to use the steepest descent curves, and the next lemma says that we can always find curves passing through $\x$ in the right directions that do the job. The main point is that during the simultaneous deformation of the $\z$ and $\e$-contours no singularity of the integrand is passed. This means that (except at $\x$) the $\e$-contour is strictly inside the $\z$-contour, 1 is between the two, and if a ray from zero hits meets the $\z$-contour at $\z$ and the $\e$-contour at $\e$, then the ratio $\e/\z$ is strictly greater than $\t$. Thus we will have to make this ratio as close to one as desired.

\noi{\bf Lemma 5}. There are disjoint closed curves $\G_\e$ and $\G_\z$ with the following properties.

\noi(i) The part of $\G_\e$ in a neighborhood $N_\e$ of $\e=\x$ is a pair of rays from $\x$ in the directions $\pm \pi/3$ and the part of $\G_\z$ in a neighborhood $N_\z$ of $\z=\x$ is a pair of rays from $\x-t^{-1/3}$ in the directions $\pm 2\pi/3$. 

\noi(ii) For some $\dl>0$ we have Re$(\ps_1(\z))<-\dl$ on $\G_\z\backslash N_\z$ and Re$(\ps_1(\e))>\dl$ on $\G_\e\backslash N_\e$.

\noi(iii) The circular $\e$ and $\z$-contours for $J$ can be simultaneously deformed to $\G_\e$ and $\G_\z$, respectively, so that during the deformation the integrand in (\ref{J}) remains analytic in all variables. 

  \begin{figure}
\begin{center}
 \includegraphics[scale=.75]{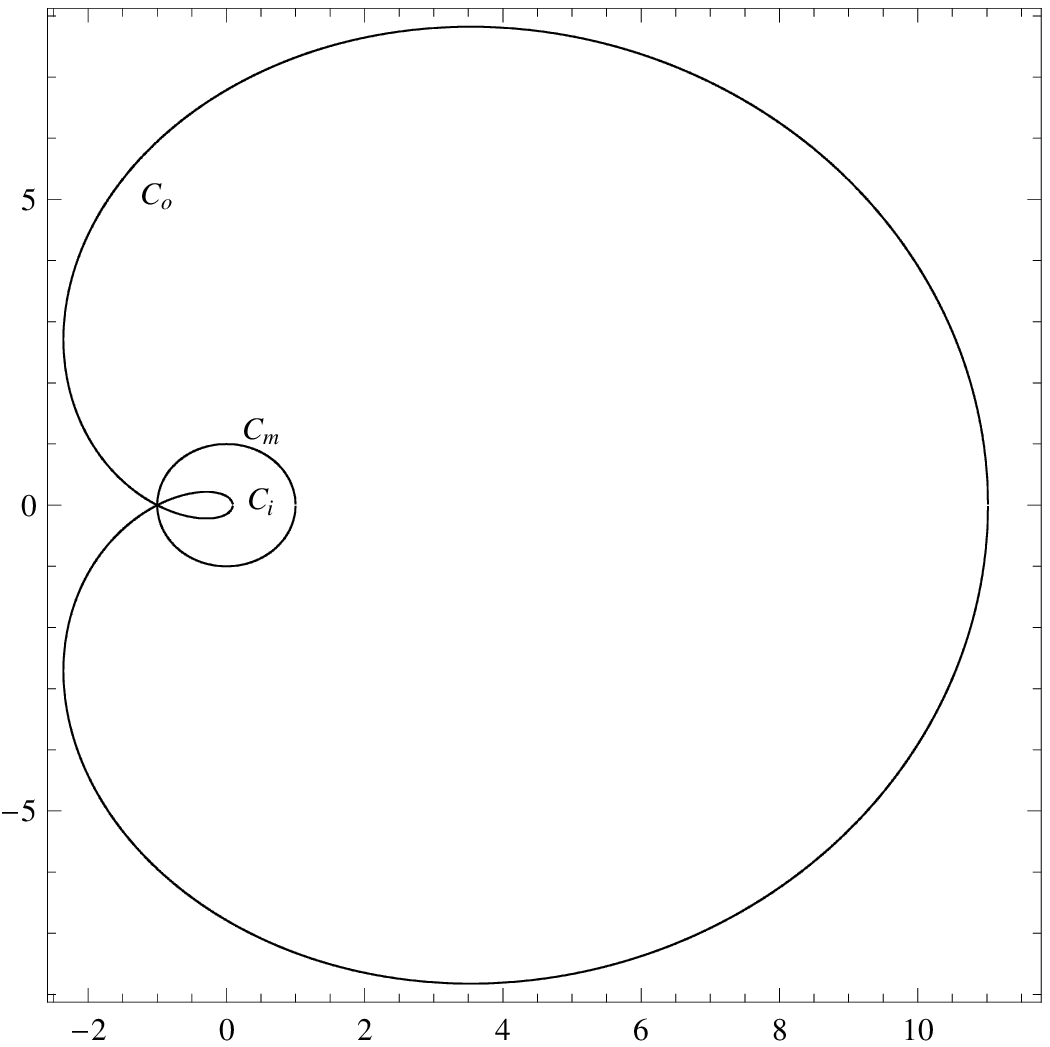}
 \caption{Curves  $C_o$, $C_m$ and $C_i$ defined by Re$(\psi_1)=0$. \label{rePsiZero}}
 \end{center}
 \end{figure} 

\noi{\bf Proof}.\footnote{The reader satisfied with an assumption that $\t$ is small enough need not read what follows.} From the local behavior of $\ps_1$ near $\x$,
\be\ps_1(\z)\sim-c_3^3\,(\z-\x)^3 /3,\label{ps1}\ee
and its global behavior we see that the set where Re$(\ps_1)=0$ consists of three closed curves meeting at $\x$. (See Fig.~\ref{rePsiZero}.) One, which we call $C_i$ since it is the inside one, has the tangent directions $\pm \pi/6$ at $\x$ and meets the real line at a point in $(0,\,1)$; another, which we call $C_m$ because it is the middle one, has the tangent directions $\pm \pi/2$ at $\x$ and meets the real line at 1; the third, which we call $C_o$ since it is the outside one, has the tangent directions $\pm 5 \pi/6$ at $\x$ and meets the real line at a point in $(1,\,\iy)$. We have Re$(\ps_1)<0$ inside $C_i$, Re$(\ps_1)>0$ between $C_i$ and $C_m$, Re$(\ps_1)>0$ between $C_m$ and $C_o$, and Re$(\ps_1)>0$ outside $C_o$. (All these may be seen by taking appropriate points in the regions and using the fact that they are connected.) Our curves $\G_\e$ and $\G_\z$ will be very close to $C_m$, the first inside it and the second outside it.

\begin{figure}
\begin{center}
 \includegraphics[scale=.75]{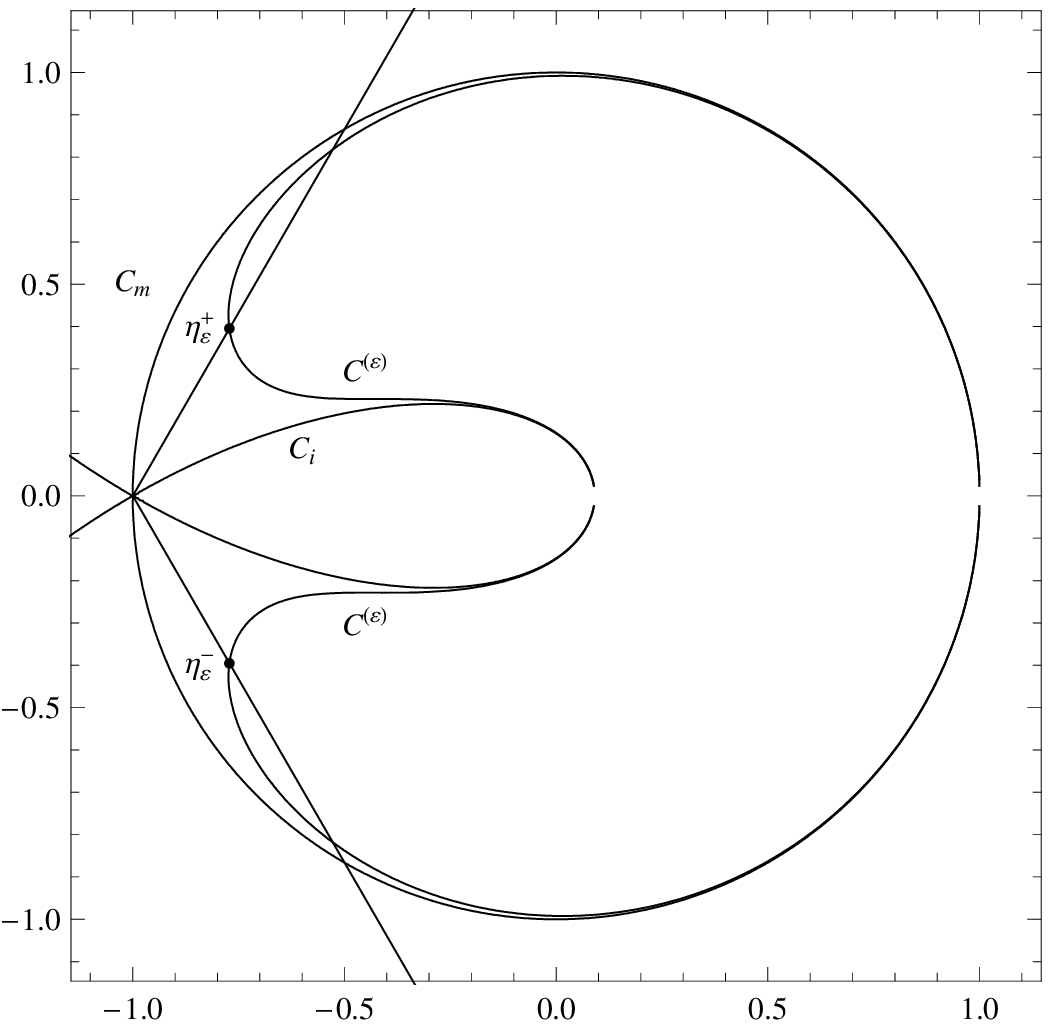}
 \caption{Curves $C^{(\ep)}$, $C_m$, $C_i$ and rays $\arg(\e-\x)=\pm \pi/3$ used in the construction of $\G_\e$. \label{rePsiPlus}}
 \end{center}
 \end{figure}

The set where Re$(\ps_1)=\ep$, with $\ep$ small and positive, consists of two curves, one lying between $C_i$ and $C_m$ and tangent to $C_m$ at $\e=1$, and the other outside $C_0$. We are interested in the first, which we call $C^{(\ep)}$. (See Fig.~\ref{rePsiPlus}.) Except for a neighborgood of $\x$, one part of $C^{(\ep)}$ is very close to $C_m$ and inside it and the other very close to $C_i$ and outside it. These are joined near $\x$ by smooth curves. The rays $\arg(\e-\x)=\pm \pi/3$ meet $C^{(\ep)}$ at points $\e_\ep^+$ and $\e_\ep^-$ close to $\x$. The curve $\G_\e$ is described as follows: it goes from $\x$ in the direction $-\pi/3$ until $\e_\ep^-$, then it takes a right turn and goes counterclockwise around $C^{(\ep)}$ (it will be very close to $C_m$ the while) until $\e_\ep^+$, and then it goes backwards along the ray with direction $\pi/3$ until returning to $\x$. (Actually, we modify this by making a semi-circular indentation around $\e=1$ to the left.)

\begin{figure}
\begin{center}
 \includegraphics[scale=.8]{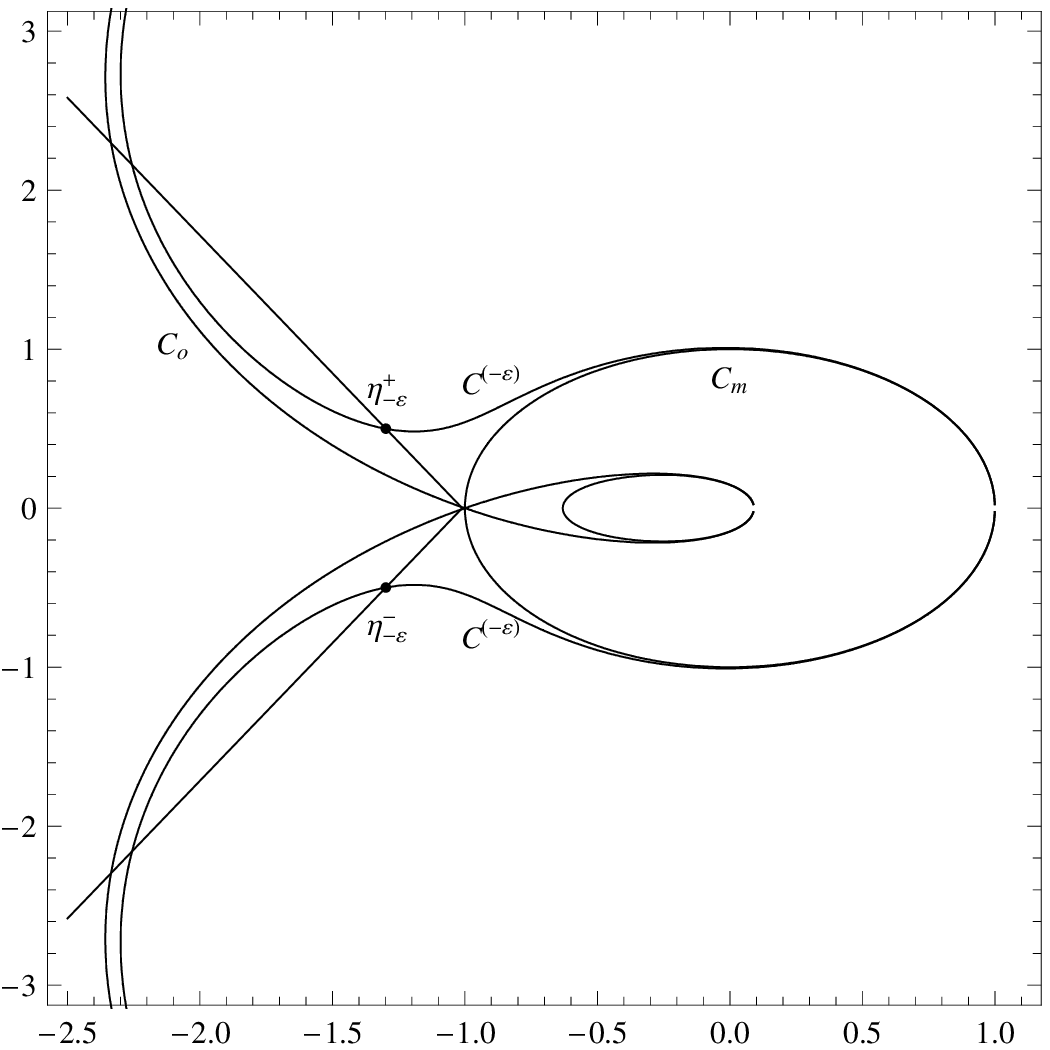}
 \caption{Curves $C^{(-\ep)}$, $C_m$, $C_o$ and rays $\arg(\z-\x+t^{-1/3})=\pm 2\pi/3$ used in the construction of $\G_\z$. \label{rePsiMinus}}
 \end{center}
 \end{figure}

The curve $\G_\z$ is obtained similarly. The set where Re$(\ps_1)=-\ep$ consists of two curves, one lying inside $C_i$ and the other between $C_m$ and $C_o$ and tangent to $C_m$ at $\z=1$. We are interested in the second, which we call $C^{(-\ep)}$. (See Fig.~\ref{rePsiMinus}.) Except for a neighborhood of $\x$, one part of $C^{(-\ep)}$ is very close to $C_m$ and outside it and the other very close to $C_o$ and inside it. These are joined near $\x$ by smooth curves. The rays $\arg(\z-\x+t^{-1/3})=\pm 2\pi/3$ meet the curves at points $\e_{-\ep}^+$ and $\e_{-\ep}^-$ near $\x$. The curve $\G_\z$ is described as follows: it goes from $\x-t^{-1/3}$ in the direction $-2\pi/3$ until $\e_{-\ep}^-$, then it takes a left turn and goes counterclockwise around $C^{(-\ep)}$ (it will be very close to $C_m$ the while) until $\e_{-\ep}^+$, and then it goes backwards along ray with direction $2\pi/3$ until returning to $\x-t^{-1/3}$. (We modify this by making a small semi-circular indentation around $\z=1$ to the right.) 

Let us see why the three stated conditions are satisfied. The first is obvious. The bounds in the second are clear on the curved parts of the contours, and is easy to see from (\ref{ps1}) on the line segments near $\x$. As for (iii), the $\z$ and $\e$-contours start out as just outside and just inside the unit circle, respectively. We may simultaneously deform these contours to just outside and inside $C_m$, respectively, without passing any singularity of the integrand in (\ref{J}), as long as the contours remain close enough to each other (and bounded away from zero). Then a further small deformation takes them to $\G_\z$ and $\G_\e$. \hfill$\Box$

\noi{\bf Proof of Theorem 3}. By part (iii) of the lemma and Proposition 1 the determinant is unchanged if $J$ acts on $\G_\e$ and the integral in (\ref{J}) is over $\G_\z$.
The operator $\m J$ is the product $AB$, where $A:L^2(\G_\z)\to L^2(\G_\e)$ and $B:L^2(\G_\e)\to L^2(\G_\z)$ have kernels
\[A(\e,\,\z)={e^{\ps(\z)}\ov \z-\e},\ \ \ 
B(\z,\,\e)={\m\,f(\m,\z/\e)\ov \e\,e^{\ps(\e)}}.\]

Aside from the factors involving $\ps$ both kernels are uniformly $O(t^{1/3})$, due to the fact that the $\z$-contour was shifted to the left by $t^{-1/3}$ near $\x$. It follows from this and (ii) that if we restrict the kernels to either $\z\in \G_\z\backslash N_\z$ or $\e\in \G_\e\backslash N_\e$ the resulting product has exponentially small trace norm. So for the limit of the determinant we may replace the contours by their portions in $N_\z$ and $N_\e$, which are rays. Using (\ref{psi}) we see that we may further restrict $\e$ and $\z$ to 
$t^{-a}$-neighborhoods of $\x$ as long as $a<1/3$, because with either variable outside such a neighborhood the product has trace norm $O(e^{-\dl\, t^{1-3a}})$.

On these segments of rays we make the replacements
\[\e\to \x+c_3\inv\,t^{-1/3}\,\e,\ \ \ \e'\to\x+c_3\inv t^{-1/3}\,\e',\ \ \ \z\to \x+c_3\inv\,t^{-1/3}\,\z.\]
The new $\e$-contour consists of the rays from 0 to $c_3\,t^{1/3-a}\,e^{\pm\pi i/3}$ while the new $\z$-contour consists of the rays from $-c_3$ to $-c_3+c_3\,t^{1/3-a}\,e^{\pm2\pi i/3}$.
In the rescaled kernels the factor $1/(\z-\e)$ in $A(\z,\,\e)$ remains the same. Because near $z=1$,
\[f(\m,\,z)=O\({1\ov|1-z|}\)\ \ {\rm and}\ \ f(\m,\,z)={\m\inv\ov 1-z}+O(1),\] 
the factor $\m\,f(\m,\z/\e)/\e$ in $B(\e,\,\z)$ becomes 
\be O\({1\ov|\e-\z|}\)\ \ {\rm and}\ \ {1\ov\e-\z}+O(t^{-1/3})\label{fest}\ee
after the rescaling. (The $\m$ and $\e$ appearing as they do is {\it very nice}.) 

As for the factors $e^{\ps(\z)}$ and $e^{-\ps(\e)}$ we see from (\ref{psi}) that for some $\dl>0$ after scaling they are $O(e^{-\dl\,|\z|^3})$ and $O(e^{-\dl\,|\e|^3})$, respectively, on their respective contours. Thus  the rescaled kernels are bounded by constants times
\[{e^{-\dl\,|\z|^3}\ov |\z-\e|},\ \ \ {e^{-\dl\,|\e|^3}\ov |\e-\z|},\]
respectively, which are Hilbert-Schmidt, i.e., $L^2$. (Notice that after the scaling $\z-\e$ becomes bounded away from zero.) It follows that convergence in Hilbert-Schmidt norm of the rescaled operators $A$ and $B$, and so trace norm convergence of their product, would be a consequence of pointwise convergence of their kernels.

The error term in (\ref{fest}) goes to zero pointwise. If also $a>1/4$, which we may assume, the error terms in (\ref{psi}) go to zero and we see that the kernels have pointwise limits 
\[{e^{-\z^3/3+s\z}\ov \z-\e},\ \ \ {e^{\e^3/3-s\e}\ov \e-\z},\]
respectively. Therefore we have found  for $\m J$ the limiting rescaled kernel
\be\int_{\G_\z} {e^{-\z^3/3+s\z+(\e')^3/3-s\e'}\ov(\z-\e)\,(\e'-\z)}\,d\z.\label{laJ}\ee
The four rays constituting the rescaled contours $\G_\z$ and $\G_\e$ in the limit go to infinity: the limiting $\G_\z$ consists of the the rays from $-c_3$ to $-c_3+\iy\,e^{\pm2\pi i/3}$ while the limiting $\G_\e$ consists of the the rays from $0$ to $\iy\,e^{\pm \pi i/3}$. 

For $\z\in\G_\z$ and $\e'\in\G_\e$ we have Re$\,(\z-\e')<0$, so we may write
\[{e^{s(\z-\e')}\ov\e'-\z}=\int_s^\iy e^{x(\z-\e')}\,dx.\]
Hence (\ref{laJ}) equals 
\[\int_s^\iy\int_{\G_\z}{e^{-\z^3/3+(\e')^3/3+x(\z-\e')}\ov\z-\e}\,d\z\,dx.\]
The operator may be written as a product $ABC$ where the factors have kernels
\[A(\e,\,\z)={e^{-\z^3/3}\ov\z-\e},\ \ \ B(\z,\,x)=e^{x\z},
\ \ \ C(x,\,\e)=e^{-x\e+\e^3/3}.\]
These are all Hilbert-Schmidt. The operator $CAB$, which has the same Fredholm determinant, acts on $L^2(s,\,\iy)$ and has kernel
\[\int_{\G_\z}\int_{\G_\e}C(x,\,\e)\,A(\e,\,\z)\,B(\z,\,y)\,d\e\,d\z\]
\[=
\int_{\G_\z}\int_{\G_\e}{e^{-\z^3/3+\e^3/3+y\z-x\e}\ov\z-\e}\,d\e\,d\z
=-\KA(x,\,y),\]
where
\[\KA(x,\,y)=\int_0^\iy \A(z+x)\,\A(z+y)\,dz.\footnote{The reason the double integral equals $-\KA(x,y)$ is that applying the operator $\partial/\partial x+\partial/\partial y$ to the two kernels
gives the same result, $\A(x)\,\A(y)$, so they differ by a function of $x-y$. Since both kernels go to zero as $x$ and $y$ go to $+\iy$ independently this function must be zero.}
\]
Hence
\[\det(I+\m\, J)\to\det\(I-\KA\,\ch_{(s,\,\iy)}\)=F_2(s).\]
The convergence is clearly uniform for $\m$ on its fixed circle, and it is easy to see that it is uniform in the neighboorhood of any fixed $\s$ and therefore for $\s$ in any compact subset of $(0,\,1)$. This completes the proof.\qed

\begin{center}{\bf Acknowledgment}\end{center}

This work was supported by the National Science Foundation through grants DMS-0553379 (first author) and DMS-0552388 (second author).


\begin{thebibliography}{99}

\bibitem{BS} Bal\'azs, M., Sepp\"al\"ainen, T.: Order of current variance and diffusivity in
the asymmetric simple exclusion process.  preprint, arXiv:0608400.

\bibitem{GK} Gohberg, I. C,  Krein, M. G.: \textit{Introduction to the Theory of Linear Nonselfadjoint Operators}. Transl. Amer. Math. Soc. {\bf 13}, Providence (1969).

\bibitem{Jo} Johansson, K.: Shape fluctuations and random matrices.
Commun.\ Math.\ Phys. \textbf{209}, 437--476  (2000).

\bibitem{KPZ} Kardar, M., Parisi, G., Zhang, Y-C.: Dynamic scaling of growing interfaces. Phys.\ Rev.\ Lett. \textbf{56}, 889--892 (1986).

\bibitem{Li} Liggett, T.M.: \textit{Interacting Particle Systems}. [Reprint of the 1985 original.] Berlin, Springer-Verlag, 2005.

\bibitem{QV} Quastel, J., Valk\'o, B.: $t^{1/3}$ superdiffusivity of finite-range asymmetric exclusion processes on Z, Commun.\ Math.\ Phys. \textbf{273}, 379--394  (2007).

\bibitem{TW1} Tracy, C.~A.,  Widom, H.: Level-spacing distributions and the Airy kernel. Commun. Math. Phys. {\bf 159}, 151--174 (1994).

\bibitem{TW} Tracy, C.~A.,  Widom, H.:  A Fredholm determinant representation in ASEP. J.~Stat. Phys. {\bf 132}, 291--300 (2008). 


\end{thebibliography}
\end{document}